\documentclass{svproc}
\usepackage{url}

\usepackage{makeidx}         
\usepackage{graphicx}        
\usepackage{amsmath, amsfonts, amssymb}
\usepackage{color}
\usepackage{bm}
\usepackage{svg}
\usepackage{multicol}        
\usepackage{booktabs}

\usepackage[colorlinks=true, linkcolor=blue, urlcolor=blue, citecolor=blue]{hyperref}

\begin{document}
\mainmatter              
\title{Reduced Basis Approach for Convection-Diffusion Equations with Non-Linear Boundary Reaction Conditions}
\titlerunning{Reduced Basis for PDEs with non-linear boundary conditions}  
%
\author{S.~Matera\inst{1}, C.~Merdon\inst{2}, D.~Runge\inst{2}}
\authorrunning{S.~Matera, C.~Merdon, D.~Runge} 
%
\tocauthor{Sebastian Matera, Christian Merdon, Daniel Runge}
\institute{Fritz Haber Institute of the Max Planck Society\\
Faradayweg 4-6, 14195 Berlin, Germany\\
\email{matera@fhi-berlin.mpg.de}\\
\and
Weierstrass Institute for Applied Analysis and Stochastics\\
Anton-Wilhelm-Amo-Stra{\ss}e 39, 10117 Berlin, Germany \\
\email{christian.merdon@wias-berlin.de, daniel.runge@wias-berlin.de}}

\maketitle              

\begin{abstract}
  This paper aims at an efficient strategy to solve drift-diffusion problems with non-linear boundary conditions as they appear, e.g., in heterogeneous catalysis.
  Since the non-linearity only involves the degrees of freedom   along (a part of) the boundary, a reduced basis ansatz is suggested that computes discrete Green's-like functions for the present drift-diffusion operator such that the   global non-linear problem reduces to a smaller non-linear  problem for a boundary method.
  The computed basis functions are completely independent of the non-linearities. 
  Thus, they can be reused for problems with the same differential operator and geometry. 
  Corresponding scenarios might be inverse problems in heterogeneous catalysis but also modeling the effect of different catalysts in the same reaction chamber.
  The strategy is explained for a mass-conservative finite volume method and
  demonstrated on a simple numerical example for catalytic CO oxidation.
\end{abstract}
{\small \textcolor{red}{This preprint has not undergone peer review or any post-submission
improvements or corrections. The Version of Record of this contribution is published in Finite {Volumes} for {Complex} {Applications} {X}—{Volume} 1, {Elliptic} and {Parabolic} {Problems},
and is available online at \url{https://doi.org/10.1007/978-3-031-40864-9_28}.}}
\section{Introduction}
Reduced basis approaches have gained popularity for the solution of partial differential equations (PDE), especially to address parameter dependencies \cite{Quarteroni2015,BennerMRA2017}. 
The idea is to employ only a few simulations with a general purpose discretization and high resolution to determine a (small) set of problem-specific basis functions (offline phase). 
This then allows for a fast solution of the PDE using the previously computed basis functions (online phase). 

We present a new reduced basis approach for linear convection-diffusion equations with highly non-linear flux boundary conditions. 
Such models are often good approximations in problems which involve surface chemistry, e.g., heterogeneous catalysis or electrochemistry \cite{MBHZGLR2015,fuhrmann_role_2011}. 
Our approach decomposes the problem into a set of linear problems to obtain the reduced basis and a non-linear problem which only depends on the degrees of freedom on the boundary. 
Being essentially a discrete representation of a type of Green's function of the linear operator, the reduced basis is independent of the non-linearity and can be reused for very different scenarios, e.g. inverse problems for parametrizing the non-linear boundary condition, which is a common task in surface chemistry applications. The connection to a type of Green's functions renders the approach closely related to boundary integral methods like \cite{trauble_nonlinear_2007}, but without explicit knowledge of the former.
In principle, the problem can be discretized by any method, e.g., finite volumes or finite element methods. 
Here, we demonstrate the approach with Voronoi finite volume methods, which allow for advantageous structural properties like conservation of mass or the non-negativity of $Y$ \cite{FUHRMANN2001201}. 
An algebraic formulation reveals that the reduced basis solution agrees with the solution of the fully coupled problem discretized with the same underlying method. 

The rest of this paper is structured as follows. 
\hyperref[sec:model_problem]{Section~}\ref{sec:model_problem}  introduces the model problem and lays out its finite volume discretization.
\hyperref[sec:reduced_basis]{Section~}\ref{sec:reduced_basis} concerns the design and implementation of the set of reduced basis functions. 
\hyperref[sec:numerics]{Section~}\ref{sec:numerics} demonstrates the approach in one model application.
Finally, \hyperref[sec:outlook]{Section~}\ref{sec:outlook} gives an outlook to further aspects and target applications.

\section{Model problem and FV discretisation}
\label{sec:model_problem}

For a given domain $\Omega$,
a (divergence-free) velocity field $\bm{v}$,
a positive diffusion coefficient $D \in \mathbb{R}$,
and a right-hand side $f \in L^2(\Omega)$, the model problem seeks
some function
$Y \in H^1(\Omega)$ such that
\begin{align}\label{eqn:model_problem}
  \mathrm{div} (Y \bm{v} - D \nabla Y) = f.
\end{align}
On the boundary $\Gamma := \partial \Omega = \Gamma_{\mathrm{in}} \cup \Gamma_{\mathrm{out}} \cup \Gamma_{\mathrm{nl}} \cup \Gamma_{\mathrm{rest}} $, various boundary conditions may
apply. Here, we assume some inlet condition
\begin{align*}
   Y = Y_\text{in} \quad \text{along} \quad \Gamma_\mathrm{in},
\end{align*}
an outflow condition
\begin{align*}
    \bm{n} \cdot (Y\bm{v}-D\nabla Y) &= Y\bm{v}\cdot \bm{n} \quad \text{along} \quad \Gamma_{\mathrm{out}},
\end{align*}
and a non-linear boundary condition
\begin{align}\label{eqn:reaction_boundary}
  \bm{n} \cdot (Y \bm{v} - D \nabla Y) = R(Y|_{\Gamma_\mathrm{nl}})
  \quad \text{along} \quad \Gamma_\mathrm{nl}
\end{align}
for some given functional $R$.
On the remaining part of the boundary $ \Gamma_{\mathrm{rest}} $ we assume homogeneous Neumann boundary
conditions for simplicity. But this approach can easily be extended to linear Robin boundary conditions.

\medskip
Here, we consider a finite volume discretization of the model problem
based on some boundary-conforming Delaunay triangulation
$\mathcal{T}$, which computes a piecewise constant approximation $Y_h \in P_0(\mathcal{K})$ to $Y$ with respect to the set of open, convex
Voronoi cells $\mathcal{K}$. Each cell $K \in \mathcal{K}$
has some associated collocation point $\bm{x}_{K} \in \overline{K} $.
The subset of cells at the boundary ${\Gamma_\mathrm{nl}}$
are denoted by $\mathcal{K}_\mathrm{nl}$ and their associated
collocation points are located on the boundary $\Gamma_\mathrm{nl}$.
Details on the implementation can be found, e.g., in the documentation of
the Julia package \texttt{VoronoiFVM.jl} \cite{fuhrmann-voronoifvm} which was also used
as a basis for the implementation of the reduced bases.

Note, that we employ an exponential fitting flux discretization that ensures desirable structural properties
like non-negativity and a maximum principle for $Y$
under certain conditions, e.g., if $R = 0$ and $\bm{v}$ is divergence-free \cite{FUHRMANN2001201}.

\section{Reduced basis approach}
\label{sec:reduced_basis}
This section describes the main idea to
speed up the computation of the solution of the model
problem with the help of a reduced basis related to the boundary
degrees of freedom of the non-linear boundary $\mathcal{K}_\mathrm{nl}$.
For this, observe that the discrete solution $Y_h$ can be decomposed into
\begin{align*}
  Y_h = Y_0 + Y_\mathrm{nl} := Y_0 + \sum_{K \in \mathcal{K}_\mathrm{nl}} \alpha_K Y_K
\end{align*}
where $Y_0$ solves the discretized linear sub-problem for $R=0$ and
each $Y_{K}$ solves the discretized fully linear problem
\begin{align}
  \mathrm{div} (Y_K \bm{v} - D \nabla Y_K) & = 0\nonumber\\ 
  \bm{n} \cdot (Y_K \bm{v} - D \nabla Y_K) & = \chi_K
  \quad \text{along} \quad \Gamma_\mathrm{nl}\label{eqn:boundary_condition_bk}
\end{align}
and homogeneous boundary conditions on the remaining boundary $ \Gamma \setminus \Gamma_{\mathrm{nl}} $.
Note, that $Y_{K}$ can be interpreted as a type of discrete
Green's function of the drift-diffusion operator for
its corressponding part $\partial K \subset \Gamma_\mathrm{nl}$
of the boundary.

To further investigate this on the algebraic
level, let $x_\mathrm{nl}$, $x_0$ and $x_K$
denote the coefficient vectors of the $Y_\mathrm{nl}$, $Y_0$ and $Y_{K}$ parts, respectively.
They are given by a solution of the linear systems of equations
\begin{align}\label{eqn:full_split_algebraic_system}
  A x_0 & = b_0\\
  A x_K & = b_K, \label{eqn:basis_algebraic_system}
\end{align}
where $A$ is the finite-volume discretization of the drift-diffusion
operator, $b_0$ encodes all linear ($Y$-independent) boundary
and right-hand side data,
and $b_K$ encodes the boundary condition \eqref{eqn:boundary_condition_bk}.
Note here, that the matrix $A$ is the same in all computations
and it is therefore straightforward to solve for $x_0$ and all $x_K$
efficiently and in parallel.
Also note, that for a constant inlet concentration and $f \equiv 0$,
it holds $Y_0 \equiv Y_\text{in}$ which allows to avoid the computation of $Y_0$ in that case.

To determine the coefficients $x_\mathrm{nl}$ for $Y_\mathrm{nl} = \sum_{K \in \mathcal{K}_\mathrm{nl}} \alpha_K Y_K$,
it is required to solve the non-linear system 
\begin{align}\label{eqn:global_system}
  A (x_0 + x_\mathrm{nl}) = b_0 + b_\mathrm{nl}\left(x_0 + x_\mathrm{nl}\right)
  \quad \Leftrightarrow \quad A x_\mathrm{nl} = b_\mathrm{nl}\left(x_0 + x_\mathrm{nl}\right)
\end{align}
where $b_\mathrm{nl}(x)$ encodes the finite volume discretisation
of the (non-linear) catalytic boundary data.
Inserting the decomposition into the reduced boundary basis
$x_\mathrm{nl} = \sum_{K \in \mathcal{K}_\mathrm{nl}} \alpha_K x_K$,
this is equivalent to seeking $\alpha_{K}$ such that
\begin{align*}
  \sum_{K \in \mathcal{K}_\mathrm{nl}} \alpha_K b_K = b_\mathrm{nl}\left(x_0 + \sum_{K \in \mathcal{K}_\mathrm{nl}} \alpha_K x_K\right).
\end{align*}
In a finite volume method that approximates the boundary integrals
in the assembly of $ b_K $ and $b_\mathrm{nl}(x)$ by a quadrature rule evaluating in the collocation point, the non-linear system to determine the coefficients $\alpha_K$
can be rewritten into
\begin{align}\label{eqn:reduced_system}
   \alpha_L = R\left(Y_0(\bm{x}_L) + \sum_{K \in \mathcal{K}_\mathrm{nl}} \alpha_K Y_K\left(\bm{x}_L\right)\right)
\quad \text{for all } L \in \mathcal{K}_\mathrm{nl}.
\end{align}
Here, $\bm{x}_L \in \Gamma_\mathrm{nl}$ denotes the collocation point
of the cell $L \in \mathcal{K}_\mathrm{nl}$.
Note, that this is a much smaller system to solve than the global system
\eqref{eqn:global_system}.

\begin{remark}
As usual in a reduced basis settting, computations can be split into
an offline and online phase. The offline phase computes the
coefficients of the linear part $x_0$ and the reduced basis functions
$Y_{K}$ resp. their coefficients $x_{K}$. The online phase solves the reduced system \eqref{eqn:reduced_system}
for a given function $R$.
\end{remark}

\begin{remark} \label{rem:compression_basis}
   The amount of work in the online phase can be further reduced
   by combining basis functions to larger ones or applying
   other compression techniques. In some applications it might
   be sufficient to make use of only one basis function that combines all
   $Y_{K}$ into a single basis function. Then \eqref{eqn:reduced_system} turns into a single equation.
\end{remark}

\begin{remark}
For solving \eqref{eqn:reduced_system} only boundary values of the reduced basis functions are needed. There is no need to store the whole vector $x_{K}$. In case volume information is needed, e.g. for plotting or evaluating quantities of interest, the full solution can be
obtained from a linear solve, where the right-hand side vector $\sum_{K \in \mathcal{K}_\text{nl}} \alpha_K b_K$ is used.
\end{remark}

\section{Model application and numerical example}
\label{sec:numerics} 

This section studies a simple, but realistic model application
that is based on the catalytic CO oxidation according to the reaction equation
\begin{align*}
  2\, \mathrm{CO} + \mathrm{O}_2 \rightarrow 2 \, \mathrm{CO}_2
\end{align*}
in a two-dimensional channel domain $\Omega := (0,5) \times (0,1)$ at a small
catalytic boundary $\Gamma_\mathrm{nl} := {0} \times (2,3)$.
The involved species mass fractions $\bm{Y} := (Y_\mathrm{CO}, Y_\mathrm{O_2},\\ Y_\mathrm{CO_2})$
with inlet mass fractions $\bm{Y}_\mathrm{in} = (0.2,0.8,0)$
at $\Gamma_\mathrm{in}$ are advected by a Hagen--Poiseuille flow
$\bm{v}(x,y) := v_\mathrm{in}(y(y-1), 0)^T$ from the inlet
$\Gamma_\mathrm{in} := \lbrace 0 \rbrace \times (0,1)$ to the outlet
$\Gamma_\mathrm{out} := \lbrace 5 \rbrace \times (0,1)$.
For the rest of the boundary $\Gamma_\mathrm{inert}$, inert wall boundary conditions are prescribed.
For simplicity, \emph{mass action kinetics} at the non-linear boundary
are assumed, which result in the reaction function
\begin{align*}
   R(Y_\mathrm{CO}, Y_\mathrm{O_2}, Y_\mathrm{CO_2}) = k\, (Y_{\mathrm{CO}})^2 (Y_{\mathrm{O}_2})^1,
\end{align*}
where $ k $ is a \emph{reaction rate constant}.
Note, that all three species are involved and their dynamics are
coupled through \emph{stoichiometric coefficients} according to the reaction above.
Altogether, we seek mass fractions $\bm{Y} $ that satisfy
\begin{align}
  \mathrm{div}(\bm{Y}\bm{v} - D \nabla \bm{Y}) & = 0 && \text{in } \Omega \nonumber \\
  \bm{Y} & = \bm{Y}_\mathrm{in} &&\text{along } \Gamma_\mathrm{in}\nonumber\\
   \bm{n} \cdot (\bm{Y} \bm{v} - D \nabla\bm{Y}) & = 0 && \text{along } \Gamma_\mathrm{inert}\nonumber\\
  \bm{n} \cdot (\bm{Y} \bm{v} - D \nabla\bm{Y}) & = \bm{Y}\bm{v}\cdot \bm{n} && \text{along } \Gamma_\mathrm{out}\nonumber\\
  \bm{n} \cdot (\bm{Y} \bm{v} - D \nabla\bm{Y}) & = R(\bm{Y}) \left[
    -2,
    -1,
    1\right]^T && \text{along } \Gamma_\mathrm{nl}.\label{eq:COOxNLBC}
\end{align}
For the sake of simplicity, we assume
that $D$ is a positive scalar, i.e.,
that the diffusion coefficients for all species
are the same and that there is no cross-diffusion.

\begin{table}[]
\setlength{\tabcolsep}{12pt}
\centering
\resizebox{\columnwidth}{!}{%
\begin{tabular}{ccc} \toprule
    {refinement level} & {global degrees of freedom ($N$)} & {reduced basis degrees of freedom} \\ \midrule
    0  & 300 & 4 \\
    1  & 1,200 & 6 \\
    2  & 4.800  & 10 \\
    3  & 19,200  & 18 \\ 
    4  & 76,800  & 34 \\ 
    5  & 307,200   & 66  \\
    6  & 1,228,800  & 130 \\ \bottomrule
\end{tabular}
}
\caption{Numbers of degrees of freedom for the global system and the reduced basis method as functions of the refinement level} \label{tab:dofs}
\end{table}
\vspace{-1em} 
The global solution is computed using \texttt{VoronoiFVM.jl} which employs a damped Newton method where the sparse linear systems are solved using a direct solver.
For the offline phase of the reduced basis method, we employ the same linear solver.
For the online phase, we exploit the linear dependence of the three non-linear boundary conditions in \eqref{eq:COOxNLBC} which reduces the number of degrees of freedom per cell to one instead of three and automatically ensures the stoichiometry. 
The non-linear system \eqref{eqn:reduced_system} is solved using the implemented Newton solver with default line search and residual norm tolerance $ \mathtt{ftol}=10^{-11} $ from the \texttt{NLsolve.jl} package \cite{mogensen_julianlsolversnlsolvejl_2020}.
We test the approach for uniformly refined meshes with $ 10 \cdot 2^{\mathrm{level}} $ nodes in each direction.
\hyperref[tab:dofs]{Table~}\ref{tab:dofs} lists the number $N$ of degrees of freedom for the global problem and the obtained reduced basis up to level $6$.

\begin{figure}
  \centering
  {
  \includegraphics[width=0.95\textwidth]{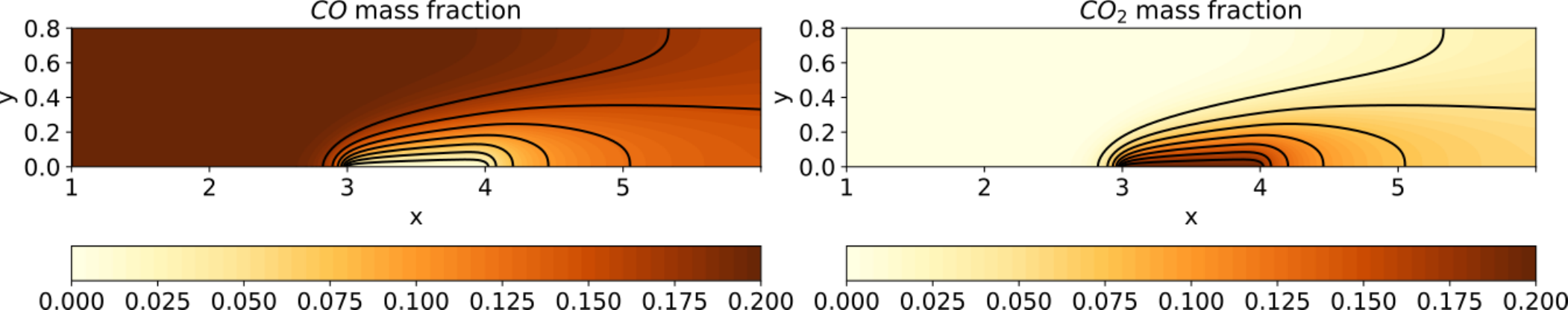}
  }
  \caption{\label{fig:massfractions}Computed mass fractions of $\mathrm{CO}$ and $\mathrm{CO}_2$ for $D = 10^{-2}$, $k = 10^{10}$ and $v_\mathrm{in} = 1$.}
\end{figure}

\hyperref[fig:massfractions]{Figure~}\ref{fig:massfractions} shows a characteristic development of the mass fractions of $\mathrm{CO}$ and $\mathrm{CO}_2$ along the catalytic surface. 
The results obtained by the reduced basis approach and the global solution agree within a tolerance close to machine precision.
While $\mathrm{CO}$ is consumed and $\mathrm{CO}_2$ is produced, a boundary layer forms whose thickness is determined by the ratio of $D$ and $v_\mathrm{in}$.
Note that the concentrations appear very uniform along the catalytic surface. This suggests that the compression described in \hyperref[rem:compression_basis]{Remark~}\ref{rem:compression_basis} could be applied here to further reduce the degrees of freedom for the online phase. 

\begin{figure}
  \centering
      \includegraphics[width = 0.45\textwidth]{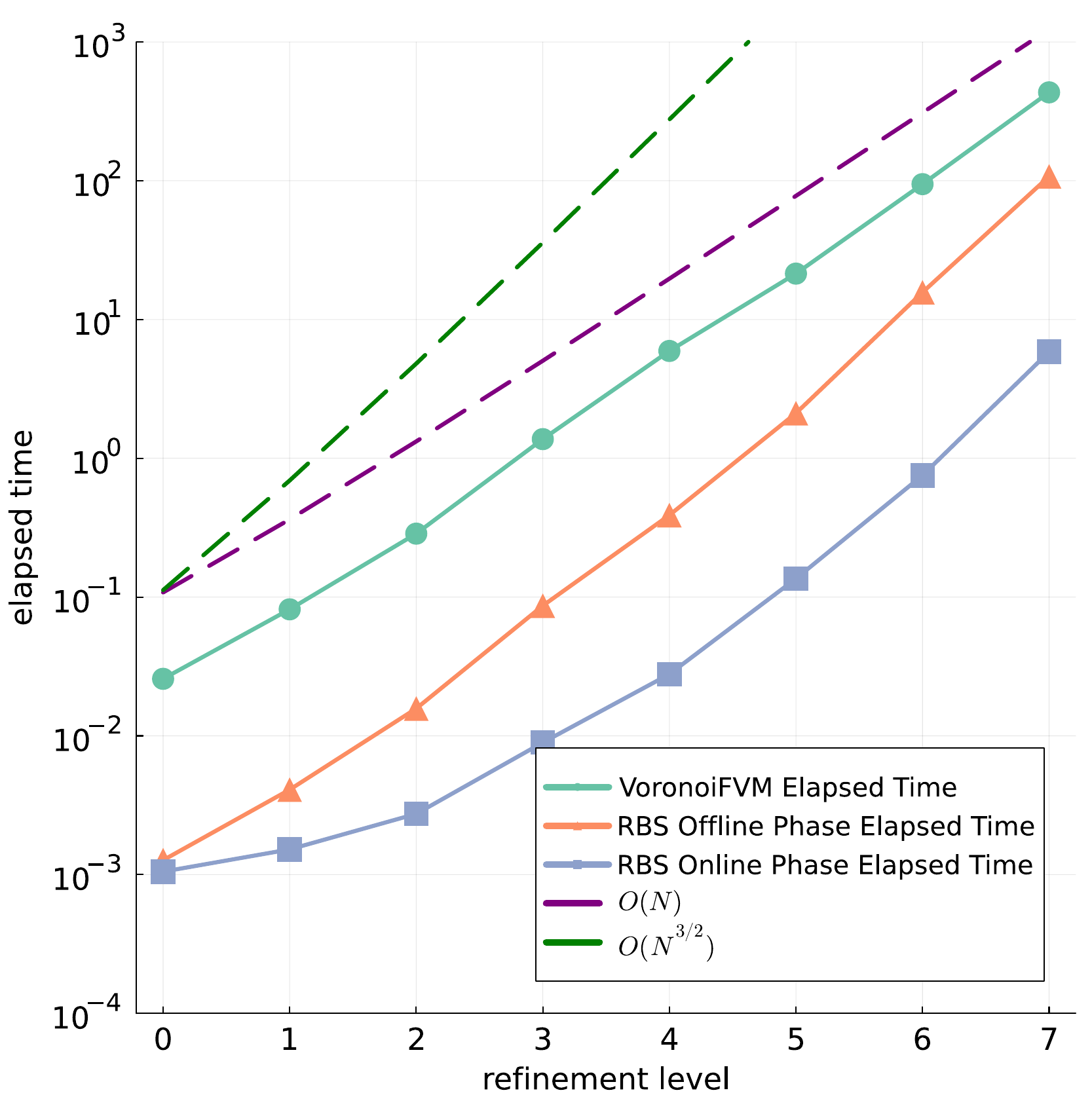}
      \hfill
      \includegraphics[width = 0.45\textwidth]{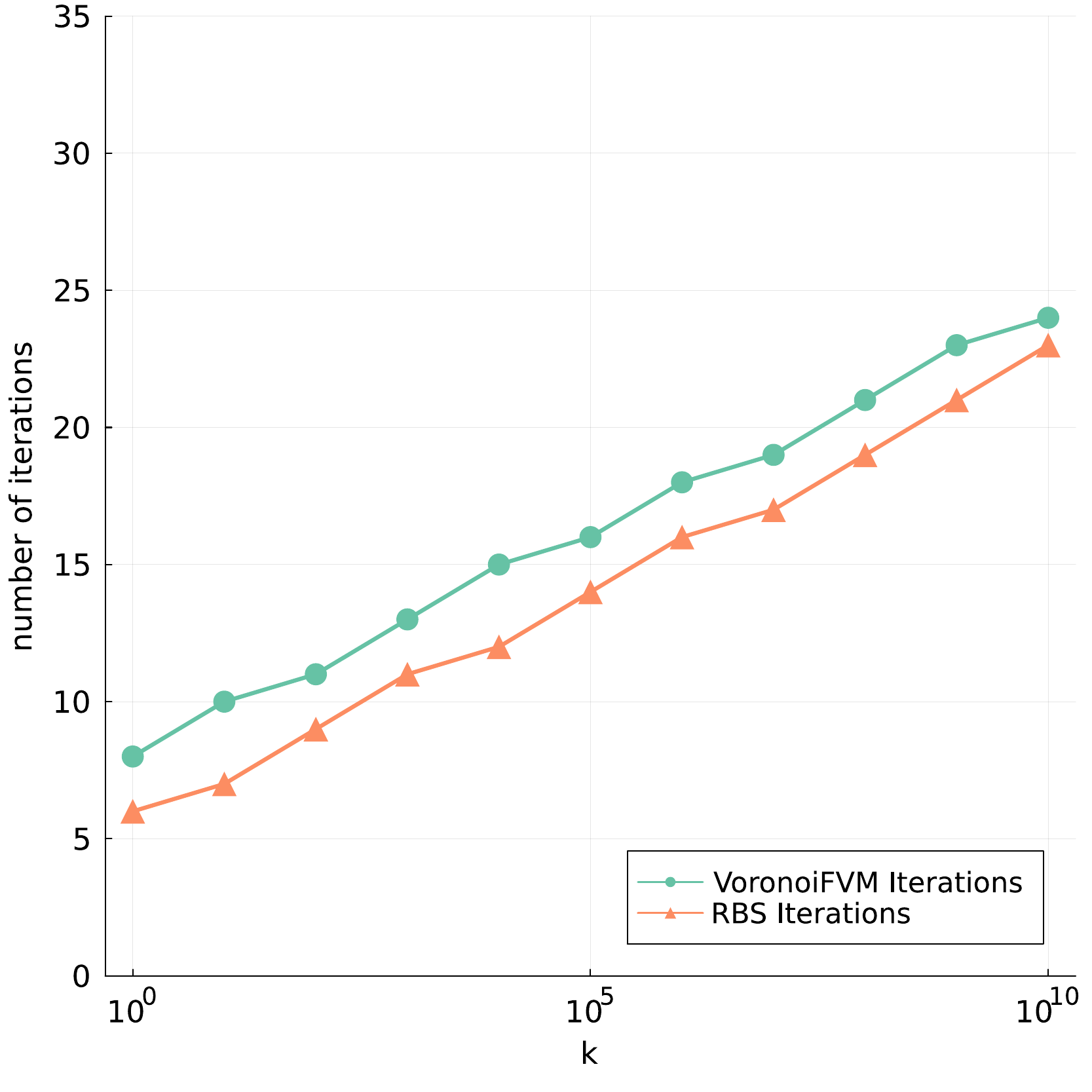}\caption{\label{fig:runtime_iterations}Runtimes as a function of the refinement level (left) for $k = 10^{10}$ and number of Newton iterations as a function of $k$ (right), both for $D = 10^{-2}$ and $v_\mathrm{in} = 1$.}
\end{figure}

\hyperref[fig:runtime_iterations]{Figure~}\ref{fig:runtime_iterations} (left) compares the runtimes of a global solve via \texttt{VoronoiFVM.jl} \cite{fuhrmann-voronoifvm} and the reduced basis scheme as functions of the number of uniform refinements for the parameters $D=10^{-2}$, $k=10^{10}$, $v_\mathrm{in} = 1$. It features the offline phase for the reduced basis setup and the online phase that solves the global problem via the reduced basis. For this and all other tested parameter settings, the online phase of the reduced basis method outperforms the global solver by about two orders of magnitude. In addition, the offline phase also comes at a significantly lower runtime.
This is partly due to the fact that the assembly of the drift-diffusion operator $A$ in \eqref{eqn:full_split_algebraic_system} is executed in every Newton step by \texttt{VoronoiFVM.jl} whereas our offline phase requires this to be executed only once.

The number of Newton iterations for a fixed $k$ in all experiments was largely
independent of the refinement level and therefore is not shown.
However, Figure~\ref{fig:runtime_iterations} (right) shows the number of Newton iterations on a fixed mesh (the finest one) versus $k$. Here, we see that the number of iterations increases with $k$, but
that the reduced basis solver always requires a comparable number of Newton iterations
as the global solver.

\section{Outlook}
\label{sec:outlook}

Here, we only discussed a simple model problem for catalytic CO oxidation and demonstrated the approach for the case where the reduced basis covers the full resolution of the global problem, simply by singling out the boundary ansatz functions.
As demonstrated, this already can reduce the computational cost dramatically, especially when many problems of the same type have to be solved on the same geometry where only the non-linearity $R(Y|_{\Gamma_\mathrm{nl}})$ is exchanged. 
The proposed methodology is thus particularly suited for inverse problems or uncertainty quantification with a parameter-dependent non-linearity. 
Since the reduced basis is completely independent of the non-linearity, this might even include qualitatively different models for the boundary reaction, e.g. for different catalyst materials. 

However, there are a number of aspects which require or allow for an extension. Particularly, this concerns the velocity field $\bm{v}$ for which no analytical expression is available in many practical applications. 
Therefore $\bm{v}$ has to be obtained from a CFD simulation. 
To ensure mass conservation within the model transport problem and also for the reduced basis, the discrete velocity field needs to be divergence-free according to the methodology derived in \cite{FUHRMANN2011530}.
 Here we plan to investigate novel, less costly divergence-free coupling strategies
 developed in the context of electrochemistry \cite{MERDON20161,FUHRMANN2019778}. 
 A way to improve efficiency is to exploit that $Y$ and thereby $R(Y|_{\Gamma_\mathrm{nl}})$ often is smooth along the boundary. This motivates to reduce the basis at the boundary, e.g.,
by combining the basis functions of neighbouring cells
or considering wavelet basis approaches. In fact, a single basis
approach has already been succesfully applied \cite{MBHZGLR2015}.

A typical class of applications is flow problems coupled with surface chemical reaction. Indeed, the current study is part of a joint effort to combine transport simulations with detailed microkinetic models  of heterogeneous catalysis.
These hybrid models shall be employed to interpret modern in situ surface characterization experiments \cite{frenken2017operando}. 
The corresponding complex employed instrumentation typically requires rather large and non-standard reaction chambers whereas the catalyst samples are rather small. 
We therefore expect the proposed methodology to be particularly effective, also because the resulting non-linear problems often are very challenging.

\section*{Acknowledgments}
The authors gratefully acknowledge the funding by the German Science Foundation (DFG) within the project “ME 4819/2-1”, the CRC 1114 “Scaling Cascades in Complex Systems” (Project No. 235221301) and under Germany’s Excellence Strategy – EXC 2008 – 390540038 – UniSysCat.

%
%
\bibliographystyle{spmpsci}
\bibliography{lit}

\end{document}